\documentclass[12pt,draftcls,onecolumn]{IEEEtran}

\IEEEoverridecommandlockouts
\usepackage{amssymb,latexsym,amsmath,epsfig,subfig,epic,
amscd,cite,oldgerm,mathrsfs,euscript,eufrak,setspace,bbm,dsfont}
\usepackage{algorithm,algpseudocode}

\newcommand{\non}{\nonumber}

\newcommand{\ba}{\begin{array}}
\newcommand{\ea}{\end{array}}
\newcommand{\be}{\begin{equation}}
\newcommand{\ee}{\end{equation}}
\newcommand{\bes}{\begin{equation*}}
\newcommand{\ees}{\end{equation*}}
\newcommand{\bea}{\begin{eqnarray}}
\newcommand{\eea}{\end{eqnarray}}
\newcommand{\beas}{\begin{eqnarray*}}
\newcommand{\eeas}{\end{eqnarray*}}
\newcommand{\ben}{\begin{enumerate}}
\newcommand{\een}{\end{enumerate}}
\newcommand{\bc}{\begin{center}}
\newcommand{\ec}{\end{center}}	
\newcommand{\bi}{\begin{itemize}}
\newcommand{\ei}{\end{itemize}}
\newcommand{\bex}{\begin{example}}
\newcommand{\eex}{\end{example}}
\newcommand{\br}{\begin{remark}}
\newcommand{\er}{\end{remark}}
\newcommand{\bd}{\begin{definition}}
\newcommand{\ed}{\end{definition}}
\newcommand{\bt}{\begin{theorem}}
\newcommand{\et}{\end{theorem}}
\newcommand{\bl}{\begin{lemma}}
\newcommand{\el}{\end{lemma}}
\newcommand{\bp}{\begin{proof}}
\newcommand{\ep}{\end{proof}}

\newcommand{\beqns}{\begin{eqnarray*}}	
\newcommand{\eeqns}{\end{eqnarray*}}		
\newcommand{\bct}{\begin{center}}		
\newcommand{\ect}{\end{center}}			
\newcommand{\benum}{\begin{enumerate}}	
\newcommand{\eenum}{\end{enumerate}}		

\newcommand{\bigblock}{\hfill \rule{1.3ex}{1.3ex}}	

\newcommand{\eps}{\epsilon}

\newcommand{\bbR}{\mathbb{R}}

\newcommand{\enma}[1]   {\ensuremath{#1}}

\newcommand{\diag}      {\enma{\mathrm{diag}}}

\newcommand{\trace}     {\enma{\mathrm{trace}}}

\newcommand{\A}{{\cal A}}

\newcommand{\eL}{{\cal L}}

\newcommand{\Pee}{{\cal P}}
\newcommand{\Q}{{\cal Q}}

\newtheorem{theorem}{Theorem}
\newtheorem{lemma}[theorem]{Lemma}

\newtheorem{proposition}[theorem]{Proposition}

\newtheorem{definition}{Definition}

\newtheorem{assumption}{Assumption}

\newtheorem{remark}{Remark}
\newtheorem{example}{Example}

\newcommand{\DefinedAs}[0]{\mathrel{\mathop:}=}

\usepackage{color}

\definecolor{maroon}{RGB}{176,48,96}
\definecolor{brown}{RGB}{165,42,42}
\definecolor{purple}{RGB}{148,0,211}
\definecolor{green}{RGB}{34,139,34}
\definecolor{orange}{RGB}{255,165,0}

\title{\bf \Large Design of Optimal Sparse Interconnection Graphs for
Synchronization of Oscillator Networks}

\author{Makan Fardad, Fu Lin, and Mihailo R. Jovanovi\'c
\thanks{Financial support from the National Science Foundation under
awards CMMI-0927509 and CMMI-0927720 and under CAREER Award CMMI-0644793
is gratefully acknowledged.}
\thanks{M.\ Fardad is with the Department of Electrical Engineering
and Computer Science, Syracuse University, NY 13244.
F.\ Lin and M.\ R.\ Jovanovi\'c are with the Department of
Electrical and Computer Engineering, University of Minnesota,
Minneapolis, MN 55455. E-mails: makan@syr.edu, fu@umn.edu,
mihailo@umn.edu.}}

\begin{document}
\maketitle

\vspace*{-1.75cm}

    \begin{abstract}
    \vspace*{-0.25cm}
We study the optimal design of a conductance network as a means for synchronizing a given set of oscillators. Synchronization is achieved when all oscillator voltages reach consensus, and performance is quantified by the mean-square deviation from the consensus value.
We formulate optimization problems that address the trade-off between synchronization performance and the number and strength of oscillator couplings. We promote the sparsity of the coupling network by penalizing the number of interconnection links. For identical oscillators, we establish convexity of the optimization problem and demonstrate that the design problem can be formulated as a semidefinite program. Finally, for special classes of oscillator networks we derive explicit analytical expressions for the optimal conductance values.
    \end{abstract}
    \vspace*{-0.25cm}
\begin{keywords}
	\vspace*{-0.25cm}
Consensus, convex relaxation, optimization, oscillator synchronization, reweighted $\ell_1$ minimization, semidefinite programming, sparse graph.
\end{keywords}

	\vspace*{-2ex}
\section{Introduction and Motivation}
\label{motive.sec}

Problems of synchronization are of interest in a variety of disciplines. In biology, examples include the synchronization of circadian pacemaker cells in the brain, pacemaker cells in the heart, and flashing fireflies and chirping crickets \cite{str00}. In engineering and applied mathematics, extensive research has been devoted to the synchronization of networks of Kuramoto oscillators and networks of power generators  \cite{dorbul12a,dorbul12b,dorjovchebulACC13,jadmotbar04,ren08}. Synchronization phenomena capture the attention of people with diverse backgrounds, as illustrated through the synchronization of mechanically coupled metronomes in the widely popular talk by Strogatz \cite{strogatz_TED}, and constitute an important part of the by now rich literature on network analysis and design \cite{easkle10,jac10,mesege10,new10}.

We consider the synchronization problem for a network of $n$ oscillators, and use the size of the conductance between any two nodes to quantify the amount of coupling between the corresponding oscillators. Our aim is to synchronize the network in a cost-effective way as far as the overall use of conductance is concerned. For oscillators subject to white-noise excitations, performance is measured using the variance of the steady-state deviation from the consensus value of oscillator voltages. We employ an ${\cal H}_2$ optimal control framework to measure the amount of synchronization and also to penalize the amount of conductance used. Additionally, in order to penalize {\em the number of} interconnection links and thus promote a sparse coupling network, we regularize the objective function with a weighted $\ell_1$ norm of the conductance matrix.

Our main contributions can be summarized as follows. We employ tools from optimal control, compressive sensing, and convex optimization to formulate the synchronization problem and design optimal sparse interconnection graphs. We develop a procedure for eliminating the marginally stable and unobservable mode which corresponds to the consensus value of oscillator voltages. Finally, we exploit problem structure to identify a class of systems for which the optimal design problem is convex, and provide a semidefinite programing formulation.

The problem of optimal controller design for large-scale and distributed systems has been considered in \cite{bampagdah02,decpag02,dandul03,lanchadan04,bamvou05,rotlal06,motjad08,farjov11}. Particular attention is paid to the problem of optimal structured control in \cite{linfarjov11a}, where the ${\cal H}_2$ norm of the closed-loop system is minimized among all controllers that respect a predetermined communication structure.
The problem of optimal sparse control is considered in \cite{linfarjov13,farlinjov11a}, where a combination of ${\cal H}_2$ norm and sparsity-promoting penalty terms is minimized with the purpose of achieving a desirable tradeoff between quadratic performance and controller sparsity. The synchronization of coupled second-order linear harmonic oscillators with local interaction is considered in \cite{ren08}. In this paper we adopt a framework which combines the optimization formulation of \cite{farlinjov11a,linfarjov13} with the oscillator network model of \cite{ren08}.

\vspace*{-2ex}
\section{Problem Formulation}
\label{state.sec}

\begin{figure}
\begin{center}
\includegraphics[scale=0.35]{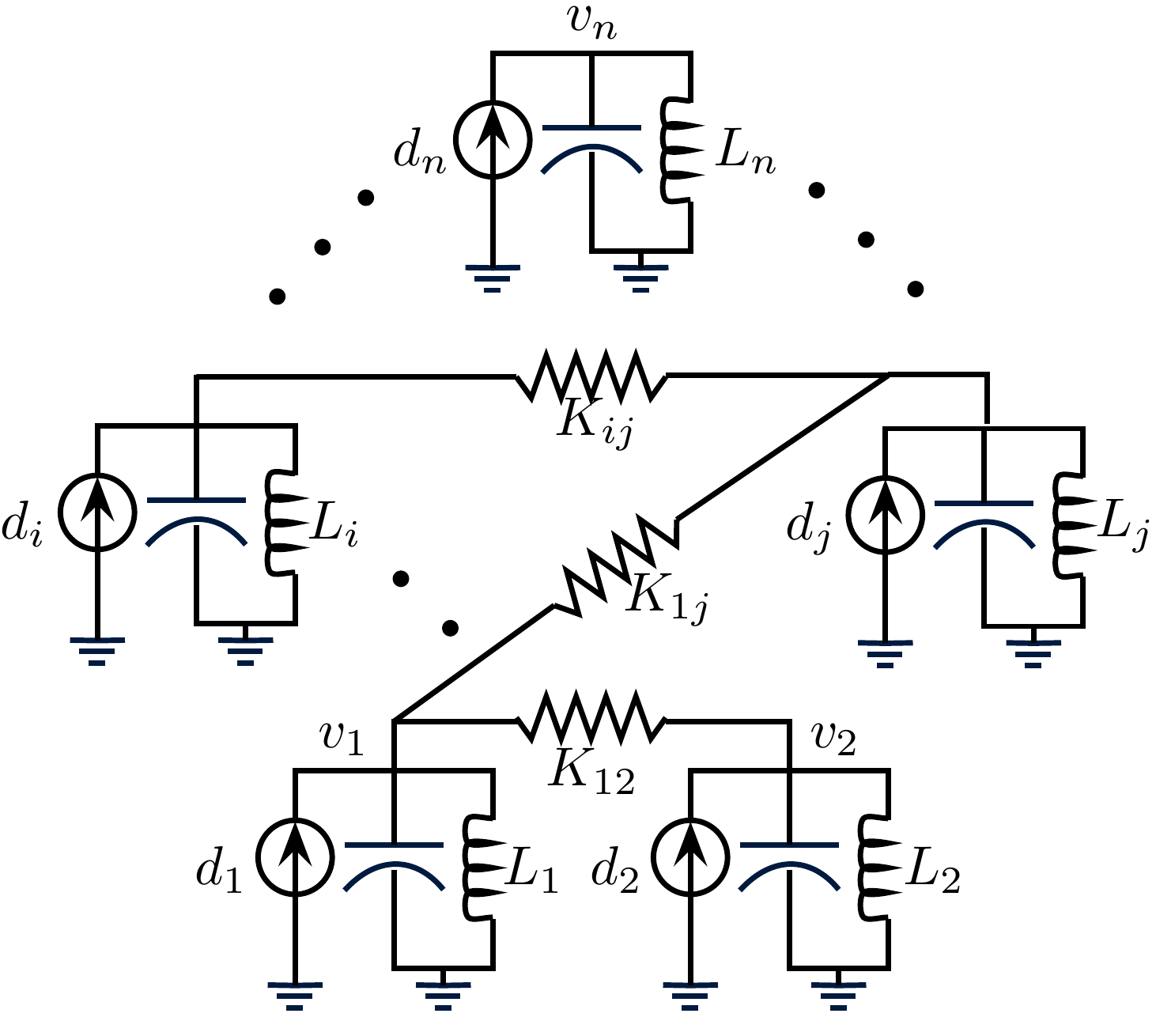}
\vspace{-0.5cm}
\end{center}
\caption{An oscillator network coupled through 
conductances described by matrix $K$.}
\label{n-oscillators.fig}
\end{figure}

Consider a network of $n$ LC-oscillators, interconnected by a set of conductances and subject to random current excitations. The conductances that connect different oscillators form the edges of an undirected (weighted) graph, with each oscillator connecting a node of the graph to the ground, as illustrated in Fig.\,\ref{n-oscillators.fig}.
For simplicity, we assume that all capacitors have unit value implying that, when considered in isolation, each oscillator resonates at frequency $\omega_i = L_i^{-1/2}$.

Let $v$ denote the column vector of node voltages. Then, taking the integral of node voltages $\int_0^t \! v$ and the node voltages $v$ as state variables, the dynamics of the entire network can be described by
\be
\dot{\psi}
\,=
\left[\ba{cc}
0 & I \\
-H & -K
\ea\right]
\psi
\,+
\left[\ba{c}
0 \\
d
\ea\right]\!,
\label{SS}
\ee
where $\psi = [\, \int_0^t \! v^T ~~ v^T \,]^T$ is the state vector, $d$ is the vector of disturbance currents injected into the nodes, 
and
\[
H \,=\, \diag \{ 1/L_i \},
~~~~~
K\! : \text{conductance matrix of node interconnections}.
\]
The conductance matrix $K$ can be thought of as a weighted Laplacian\cite{ghoboysab08}, which by default satisfies
\[
K \succeq 0,
~~~~~
K \mathds{1} = 0.
\]
Here, $\succeq$ denotes inequality in the matrix semidefinite sense and $\mathds{1}$ is the column vector of all ones. We assume that the system's graph is connected, which implies the positive definiteness of the matrix $K$ when it is restricted to the subspace $\mathds{1}^\perp$.
Concisely, we write $K \in \eL$, where
\begin{align}
\eL
&\DefinedAs
\{ 
K 
~|~
K \, = \, K^T,
~
K \, \mathds{1} \, = \, 0,
~
K \, + \, \mathds{1} \mathds{1}^T/n \, \succ \, 0,
~
K_{ij} \, \leq \, 0 ~\text{for}~ i \neq j
\}.
\label{L}
\end{align}

It is desired to find an `optimal' (in a sense to be made precise in what follows) matrix $K$ such that:
(i)
the difference in node voltages $|v_i - v_j|$ is kept small for every $i$ and $j$;
(ii)
the total amount of conductance used to connect nodes is kept small; and
(iii)
the number of links between nodes is kept small.
Objective\;(i) attempts to synchronize the oscillators by keeping the node voltages close to each other.
Objective\;(ii) tries to maintain a small level of coupling between the nodes.
Objective\;(iii) aims to obtain a sparse interconnection topology.
We note that objective\;(iii) is sometimes relaxed in this paper, for example, when a particular interconnection topology is determined {\em a priori} and optimal values of conductances are sought within that topology.

In order to place the problem of designing $K$ in the framework of optimal control theory, we rewrite system (\ref{SS}) in state-space form \cite{zhou} as
\begin{align*}
\dot{\psi}
\; &= \; 
A \, \psi \,+\, B \, d \,+\, B \, u,
\non 
\\[-0.1cm]
z 
\; &=\; 
\! \left[\ba{c}
Q^{1/2} \\
0
\ea\right] \psi 
\;+\, 
\left[\ba{c}
0 \\
R^{1/2}
\ea\right] u,  
\\[-0.1cm]
y \; &=\; C \, \psi,
\end{align*}
with ~
$
u = -K \, y
$,
$
K \in \eL
$,
and
\begin{align*}
&
A  
\; = \,
\left[\ba{cc}
0 & I \\
-H & 0
\ea\right]\!,
~~~
B 
\; = \,
\left[\ba{c}
0 \\
I
\ea\right]\!,
~~~
C
\; = \,
\left[\ba{cc}
0 & I
\ea\right]\!.
\end{align*}
The variables $d$ and $u$ respectively represent the exogenous and control inputs that enter the nodes as currents, and the variables $z$ and $y$ respectively represent the performance and measured outputs. The positive semidefinite matrix $Q$ and the positive definite matrix $R$  respectively quantify state and control weights. In this control-theoretic framework the matrix $K$ denotes the static feedback gain, which is subject to the structural constraint of being in the set $\eL$ defined in (\ref{L}). Upon closing the loop, the above problem can equivalently be written as
\begin{align}
\dot{\psi}
\,&=\,
(A - B K C) \, \psi \,+\, B \, d,
\label{CL}
\\
z
\,&=\,
\left[
\ba{c}
\!\! Q^{1/2} \!\! \\
\!\! -R^{1/2} K C \!\!
\ea
\right] \psi.
\non
\end{align}
From the above definitions of the matrices $A$, $B$, and $C$ it is easy to see the equivalence between the equations (\ref{CL}) and (\ref{SS}).

We further assume that
$
R = r \,I
$,
$
r > 0
$,
and
\be
Q
\,=
\left[\ba{cc}
0 & 0 \\
0 & Q_2
\ea\right]\!,
\label{Q}
\ee
where $Q_2$ satisfies $Q_2 \mathds{1} = 0$ and is a positive definite matrix when restricted to the subspace $\mathds{1}^\perp$, 
\[
Q_2 \mathds{1} \,=\,0,
~~
\zeta^T Q_2 \, \zeta \,>\, 0
~\text{for all}~ \zeta \neq 0 ~\text{such that}~ \zeta^T \mathds{1} = 0.
\]
To justify the structural assumptions on $Q$, we note that in order to achieve synchronization we are interested in making weighted sums of terms of the form $(v_i - v_j)^2$ small. Owing to the choice of state variables $[\, \int \! v^T ~~ v^T \,]^T$, such an objective corresponds to $Q$ matrices with the zero structure displayed in (\ref{Q}) and $Q_2$ matrices that are positive semidefinite and satisfy $Q_2 \mathds{1} = 0$. For example, in a system of two oscillators, if it is desired to make $(v_1 - v_2)^2$ small then $Q$ has the structure shown in (\ref{Q}) with
\[
Q_2
\,=
\left[\ba{c}
1 \\
-1
\ea\right]\!
\left[\ba{cc}
1 & -1
\ea\right]
=
\left[\ba{rr}
1 & -1 \\
-1 & 1
\ea\right]\!.
\]

We now state the main optimization problem addressed in this paper, and then elaborate on the details of its formulation. Consider the problem 
\be
\ba{ll}
\text{minimize}
& J_\gamma 
\,\DefinedAs\, 
\trace \,  ( P B B^T ) \,+\, \gamma \, \| W \circ K \|_{\ell_1}
\\[0.15cm]
\text{subject to}
&
(A \, - \, B K C)^T P \,+\, P (A \, - \, B K C)
\;=\;
-(Q \, + \, C^T K^T R K C)
\\[0.15cm]
&
K \, \in \, \eL,
~~~
P \, \succeq \, 0,
\ea
\label{OPT-general}
\ee
where $K$ and $P$ are the optimization variables,
$\| K \|_{\ell_1} = \sum_{i,j} |K_{ij}|$ is the $\ell_1$-norm of $K$,
$W$ is a weighting matrix,
$\circ$ denotes elementwise matrix multiplication, and $\eL$ is the set of weighted Laplacian matrices corresponding to connected graphs, as defined in (\ref{L}).

We next elaborate on the formulation of the optimization problem (\ref{OPT-general}). 
When $\gamma = 0$, the objective function
\[
J
\,\DefinedAs\,
\trace \,  ( P B B^T )
\]
determines the ${\cal H}_2$ norm, from input $d$ to output $z$, of the closed-loop system (\ref{CL}) \cite{zhou}. In a stochastic setting, the ${\cal H}_2$ norm quantifies variance amplification from $d$ to $z$ in statistical steady-state; in a deterministic setting, it quantifies the $L_2$-norm of the impulse response.
Solving (\ref{OPT-general}) for $\gamma = 0$ is closely related to the design of structured feedback gains \cite{linfarjov11a}. The condition $K \in \eL$ ensures that $K$ is a legitimate conductance matrix. 
It is important to note that due to our particular choice of the performance output $z$, by minimizing $J$ we are effectively achieving the first two of our optimal synchronization objectives outlined earlier. 
Furthermore, it has been demonstrated recently that $\ell_1$ optimization can be effectively employed as a proxy for cardinality minimization~\cite{canromtao06,canwakboy08}, where the cardinality $\mathrm{card}(K)$ of a matrix $K$ is defined as the number of its nonzero entries. Indeed, the term $\gamma \, \| W \circ K \|_{\ell_1}$ in the objective function of (\ref{OPT-general}) attempts to approximate $\gamma \, \mathrm{card}(K)$ in penalizing the number of nonzero elements of $K$, which in terms of the synchronization problem can be interpreted as penalizing the number of interconnection links. The weighting matrix $W$ can be updated via an iterative algorithm in order to make the weighted $\ell_1$ norm $\| W \circ K \|_{\ell_1}$ a better approximation of $\mathrm{card}(K)$~\cite{canwakboy08,farlinjov11a}. We next describe one such algorithm.

	\vspace*{-2ex}
\subsection{Sparsity--Promoting Reweighted $\ell_1$ Algorithm}
\label{reweight.sec}

Reference \cite{canwakboy08} introduces the {\em reweighted $\ell_1$ minimization} algorithm as a relaxation for cardinality minimization. This methodology was recently used in \cite{linfarjov13,farlinjov11a} to find optimal sparse controllers for interconnected systems. We now state the reweighted $\ell_1$ algorithm for the optimal sparse synchronization problem.\\

\begin{algorithm}
\caption{Reweighted $\ell_1$ algorithm}
\begin{algorithmic}[1]
\State	{\bf given} $\delta > 0$ and $\eps > 0$.
\For{$\mu = 1,2,\ldots$}
\State	If $\mu = 1$, set $K^\mathrm{prev} = 0$, set $W_{ij} = 1$, form $W$.
\State	If $\mu > 1$, set $K^\mathrm{prev}$ equal to optimal $K$ from previous

\hspace{-0.3cm} step, set
$
W_{ij} = 1/(|K^\mathrm{prev}_{ij}| + \delta)
$,
form $W$.
\State	Solve (\ref{OPT-general}) to obtain $K^*$.
\State	If $\| K^* - K^\mathrm{prev} \| < \eps$, {\bf quit}.
\EndFor
\end{algorithmic}
\end{algorithm}

\noindent
Henceforth in this paper, unless stated otherwise, we will only address solving the optimization problem (\ref{OPT-general}) for a given weighting matrix $W$, which corresponds to Step\,5 of Algorithm\,1.

	\vspace*{-2ex}
\subsection{Simplification of Problem {\em (\ref{OPT-general})}}
\label{simplify.sec}

We begin by exploiting the structure of the Lyapunov equation that appears in the optimization problem (\ref{OPT-general}),
\[
(A \, - \, B K C)^T P \,+\, P (A \, - \, B K C)
\;=\;
-(Q \, + \, C^T K^T R K C).
\]
Substituting the expressions for $A$, $B$, $C$, $Q$, and
\[
P
\,=
\left[\ba{cc}
P_1 & P_0 \\
P_0^T & P_2
\ea\right] \succeq\, 0,
\]
and
rewriting the equation in terms of its components gives
\begin{align}
H P_0^T \,+\, P_0 H
\,&=\,
0
\non
\\
P_0 K \,-\, P_1 \,+\, H P_2
\,&=\,
0
\label{LYP}
\\
K P_2 \,+\, P_2 K \,-\, P_0 \,-\, P_0^T
\,&=\,
Q_2 \,+\, r K^2.
\non
\end{align}
The condition $P \succeq 0$ implies that $P_1 \succeq 0$ and $P_2 \succeq 0$. 
Finally, we use the block decomposition of $P$ to simplify the objective function in (\ref{OPT-general}),
\be
\trace \, (P B B^T) \,=\, \trace \, (P_2).
\label{TP}
\ee

	\vspace*{-2ex}
\section{Case of Uniform Inductances: A Convex Problem}
\label{sdp.sec}

For networks in which all inductors have the same value, we show in this section that the optimization problem (\ref{OPT-general}) can be formulated as a semidefinite program. 

\noindent
\begin{assumption}
\label{uniform.ass}
Let all inductors have the same value, i.e.,
\be
L_i \,=\, L,
~~~
i = 1,\ldots,n,
~~~~~
H \,=\, (1/L) I.
\label{UI}
\ee
for some $L > 0$. We hereafter refer to this as the `uniform inductance' assumption.
\end{assumption}

\br
This assumption is restrictive in that all oscillators now have the same resonance frequency $\omega = L^{-1/2}$ (recall that all capacitor values are equal to one). However, the synchronization problem is still meaningful, as it forces the oscillators to reach consensus on their amplitudes and phases and oscillate in unison. It can be shown \cite{sync_arxiv} that the uniform inductance scenario provides a valuable design platform for the more general case in which different inductor values constitute small deviations from some nominal value $L_0$.
\er

From the uniform inductance assumption (\ref{UI}) it follows that $H = (1/L) I$ commutes with any matrix and therefore the first equation in (\ref{LYP}) becomes
\[
P_0 \,+\, P_0^T \,=\, 0.
\]
Hence the last equation in (\ref{LYP}) simplifies to
\be
K P_2 \,+\, P_2 K
\,=\,
Q_2 \,+\, r K^2,
\label{P2}
\ee
with $P_2 \succeq 0$. Furthermore, from (\ref{TP}) it follows that the objective in (\ref{OPT-general}) is equal to $\trace \, (P_2)$ and is independent of $P_0$ and $P_1$. The optimization problem (\ref{OPT-general}) can thus be rewritten as
\be
\ba{ll}
\text{minimize}
&
\trace \, (P_2) \,+\, \gamma \, \| W \circ K \|_{\ell_1}
\\[0.15cm]
\text{subject to}
&
K P_2 \,+\, P_2 K
\,=\,
Q_2 \,+\, r K^2
\\[0.15cm]
&
K \in \eL,
~~~
P_2 \succeq 0.
\ea
\label{OPT-uniform}
\ee
It is worth noting the close correspondence between the optimization problem (\ref{OPT-uniform}) and a related optimal sparse design problem for a network of single-integrators \cite{farlinjov11a}.

To simplify the optimization problem further, we state the following useful lemma.

\begin{lemma}
\label{zero-mode.thm}
Let $\A$ and $\Q$ be given symmetric matrices that satisfy $\A \mathds{1} = \Q \mathds{1} = 0$, and suppose that $\A$ is negative definite when restricted to the subspace $\mathds{1}^\perp$. 
For the Lyapunov equation
\be
\A^T \Pee \,+\, \Pee \A \,=\, -\Q,
\label{lyap-curly}
\ee
the following statements hold.
\bi
\item[(i)]
If $\Pee$ is a solution of the Lyapunov equation (\ref{lyap-curly}) then so is $\Pee + \alpha \, \mathds{1}  \mathds{1}^T$ for any $\alpha \in \bbR$.
\item[(ii)]
If $\Pee$ is a solution of the Lyapunov equation (\ref{lyap-curly}) and $\Q$ is positive semidefinite on $\mathds{1}^\perp$, then $\mathds{1}$ is an eigenvector of $\Pee$ and $\Pee$ is positive semidefinite on $\mathds{1}^\perp$. Furthermore, among all $\Pee \succeq 0$ that satisfy (\ref{lyap-curly}) the one with the minimum trace satisfies $\Pee \mathds{1} = 0$.
\item[(iii)]
Any solution $\Pee$ of the Lyapunov equation (\ref{lyap-curly}) satisfies
\begin{align*}
\trace \, (\Pee)
\,&=\,
p \, - \, (1/2) \, \trace \, (\Q \A^\dagger),
~~
\Pee
\,=\,
(p/n) \, \mathds{1} \mathds{1}^T \,+\, \Pee^\perp,
\end{align*}
for some $p \in \bbR$ and matrix $\Pee^\perp$ with $\Pee^\perp \mathds{1} = 0$, where $p$ is independent of $\A$ and $\Q$, and $\A^\dagger$ denotes the Moore--Penrose pseudoinverse of $\A$.  
Additionally, if $\Q$ is positive semidefinite on $\mathds{1}^\perp$ then so is $\Pee^\perp$.
\item[(iv)]
The identity
$\trace \, (\Q \A^\dagger)
=
\trace \, (\Q (\A - \mathds{1} \mathds{1}^T/n)^{-1})
$
holds and any solution $\Pee$ of the Lyapunov equation (\ref{lyap-curly}) satisfies
\begin{align*}
\trace \, (\Pee)
\,&=\,
p - (1/2) \, \trace \, (\Q (\A - \mathds{1} \mathds{1}^T/n)^{-1}),
\end{align*}
where $p \in \bbR$ is independent of $\A$ and $\Q$.
\ei
\end{lemma}

{\em Proof:}~
The proof uses a special similarity transformation to eliminate the common zero mode of $\A$ and $\Q$ from the Lyapunov equation $\A^T \Pee + \Pee \A = -\Q$; see Appendix for details.
\bigblock


\br
An important consequence of Lemma\,\ref{zero-mode.thm} is that the new description of $\trace \, (\Pee)$,
\begin{align*}
\trace \, (\Pee)
\,&=\,
p - (1/2) \, \trace \, (\Q^{1/2} (\A - \mathds{1} \mathds{1}^T/n)^{-1} \Q^{1/2}),
\end{align*}
lends itself to the application of semidefinite programing (SDP) methods, as we  demonstrate below. This is reminiscent of the results in \cite{ghoboysab08}.
\er


Applying Lemma\,\ref{zero-mode.thm} with $\A = -K$ and $\Q = Q_2 + r K^2$ to the Lyapunov equation (\ref{P2}) with $P_2 \succeq 0$ gives
\begin{align}
J
&=\,
\trace \, (P_2)
\non
\\
&=\,
p + (1/2) \, \trace \, ((Q_2 + r K^2) K^\dagger)  
\non
\\
&=\,
p + (1/2) \, \trace \, (Q_2 (K + \mathds{1} \mathds{1}^T/n)^{-1}
+ r K (I - \mathds{1} \mathds{1}^T/n))  
\non
\\
&=\,
p + (1/2) \, \trace \, (Q_2^{1/2} (K + \mathds{1} \mathds{1}^T/n)^{-1} Q_2^{1/2} + r K),
\label{TP2}
\end{align}
where $p \in \bbR$ is independent of $K$ and $Q_2$.
The details of the simplifications in (\ref{TP2}) are as follows: Since $K$ is the Laplacian of a connected graph then $K \mathds{1} = 0$ and $K$ is positive definite on $\mathds{1}^\perp$. Also, by assumption $Q_2 \mathds{1} = 0$ and $Q_2$ is positive definite on $\mathds{1}^\perp$. Thus $(Q_2 + r K^2) \mathds{1} = 0$ and $Q_2 + r K^2$ is positive definite on $\mathds{1}^\perp$. Therefore Lemma\,\ref{zero-mode.thm} applies and the first equation follows. In the second equation, $Q \mathds{1} = 0$ and the identities
$K^\dagger = (K + \mathds{1} \mathds{1}^T/n)^{-1} - \mathds{1} \mathds{1}^T/n$,
$K^\dagger K = I - \mathds{1} \mathds{1}^T/n$
are invoked. Finally, the last equation follows from $K \mathds{1} = 0$ and the trace identity $\trace \, (M_1 M_2) = \trace \, (M_2 M_1)$.

In summary, problem (\ref{OPT-general}), which has been simplified to (\ref{OPT-uniform}) using the uniform inductance assumption, is further simplified with the help of Lemma\,\ref{zero-mode.thm} and (\ref{TP2}) to obtain the equivalent problem
\be
\!\!\!\!
\ba{ll}
\text{minimize}
&\! (1/2) \, \trace \, (Q_2^{1/2} (K + \mathds{1} \mathds{1}^T/n)^{-1} Q_2^{1/2} + r K)
\, + \, \gamma \, \| W \circ K \|_{\ell_1}
\\[0.15cm]
\text{subject to}
&\!
K \, \in \, \eL,
\ea
\label{CVX}
\ee
where the parameter $p$ has been dropped from the objective, as it has no effect on the solution of the optimization problem.

	\vspace*{-2ex}
\subsection{SDP Formulation}
\label{sdp.sec}


\begin{proposition}
The optimization problem (\ref{OPT-general}), under the uniform inductance assumption, is equivalent to the semidefinite program
\be
\!\!\!\!\!\!\!\!
\ba{ll}
\mathrm{minimize}
&\!\!
(1/2) \, \trace \, (X \, + \, r K) \, + \, \gamma \, \trace \, (\mathds{1} \mathds{1}^T Y)
\\[0.15cm]
\mathrm{subject ~ to}
&\!\!
\left[\ba{cc}
X & Q_2^{1/2} \\
Q_2^{1/2} & K + \mathds{1} \mathds{1}^T/n
\ea\right] \succeq 0
\\[0.4cm]
&\!\!
M \circ K \leq 0,
~\,
K \mathds{1} = 0,
~
-Y \leq W \circ K \leq Y,
\ea
\label{SDP}
\ee
where the optimization variables are the symmetric matrices $K$, $X$ and the elementwise-nonnegative matrix $Y$, $\leq$ denotes elementwise inequality of matrices, and 
$ M \DefinedAs \mathds{1} \mathds{1}^T - I$.
\end{proposition}

{\em Proof:}~
See Appendix.
\bigblock

\noindent
We note that the optimal conductance matrix is independent of the inductance matrix $H$ when all inductances have the same value. In other words, the optimal $K$ does not depend on the oscillator parameters when all oscillators are identical.

	\vspace*{-2ex}
\subsection{Optimality Conditions for $\gamma = 0$}
\label{all2all-1.sec}

\begin{proposition}
Consider $J(K) = \trace \,  ( P(K) B B^T )$ subject to the constraints in (\ref{OPT-general}) and the uniform inductance assumption. 
Then
\begin{align*}
\nabla_{\!K} J
\,&=\, -(1/2) \, (K + \mathds{1} \mathds{1}^T/n)^{-1} Q_2 (K + \mathds{1} \mathds{1}^T/n)^{-1}
\,+\,
(r/2) \, (I - \mathds{1} \mathds{1}^T/n)
\end{align*}
In particular, setting $\nabla_{\!K} J = 0$ gives
$
K = Q_2^{1/2}/r^{1/2}
$
as a necessary and sufficient condition for the optimality of $K$. 
\end{proposition}

{\em Proof:}~
The proof follows from a straight forward application of variational methods to the expression in (\ref{TP2}) and noting that $K$ is restricted to the set $\eL$. We omit the details due to space limitations.\bigblock ~

	\vspace*{-2ex}
\section{Illustrative Examples}
\label{example.sec}


\begin{example}[Uniform All-To-All Coupling]

A problem of particular interest in oscillator synchronization is that of uniform all-to-all coupling \cite{str00}. This structured (nonsparse) problem can be easily addressed using the framework developed above. In this case every oscillator is connected to all other oscillators and all couplings have the same magnitude. This implies a particular structure on $K$, namely
\[
K \,=\,  k \, (I - \mathds{1} \mathds{1}^T/n),
~~~~~
k > 0.
\]
It is easy to see that this $K$ belongs to $\eL$.

Since the structure of $K$ is already determined, the sparsity-promoting term $\gamma \, \| W \circ K \|_{\ell_1}$ can be dropped from the objective of (\ref{CVX}), and the problem simplifies to finding the value of  $k$ that minimizes $J = \trace \, (P_2)$. 
We have 
$
K^\dagger =  (1/k) \, (I - \mathds{1} \mathds{1}^T/n)
$,
$
(K + \mathds{1} \mathds{1}^T/n)^{-1} = (1/k) \, I + (1 - 1/k) \, \mathds{1} \mathds{1}^T/n
$, 
and
\begin{align*}
2 J
\,&=\,
\trace \, (Q_2^{1/2} (K + \mathds{1} \mathds{1}^T/n)^{-1} Q_2^{1/2} + r K)
\,=\,
(1/k) \, \trace \, (Q_2)
\,+\,
r k \, (n - 1).
\end{align*}
Setting $\partial J / \partial k = 0$ we obtain
\[
k \; = \, \left( \dfrac{\trace \, (Q_2)}{(n-1) r} \right)^{1/2},
\]
thereby implying that the optimal conductance matrix is given by
\be
K 
\; = \, 
\left(
\dfrac{\trace \, (Q_2)}{(n-1) r} \right)^{1/2}
(I \, - \, \mathds{1} \mathds{1}^T/n).
\label{KA}
\ee
Notice that because of the particular structure enforced on $K$, the optimal conductance matrix depends only on the trace of $Q_2$ and not on its exact structure or its individual entries.
\end{example}


\begin{example}

In this example we consider $n = 7$ identical oscillators and design a sparse conductance matrix using the sparsity-promoting algorithm of Section~\ref{state.sec}, with the optimization problem in Step\,5 of Algorithm\,1 being (\ref{SDP}).

Let $r = 1$ and $Q_2$ denote the $7$-by-$7$ version of the matrix
\[
Q_2
\; \sim \,
\left[\ba{rrrrrrr}
\!\!  1  \!&\!  -1  \!&\!   0  \!&\!   0  \!\!\\
\!\! -1  \!&\!   2  \!&\!  -1  \!&\!   0  \!\!\\
\!\!  0  \!&\!  -1  \!&\!   2  \!&\!  -1 \!\!\\
\!\!  0  \!&\!   0  \!&\!  -1  \!&\!   1 \!\!
\ea\right].
\]
The optimal conductance matrices $K_\gamma$, for different values of $\gamma$, are given below. For all computations we used \texttt{CVX}, a package for specifying and solving convex programs \cite{CVX,gb08}. As expected, for $\gamma = 0$ we recover $K_0 = Q_2^{1/2}/r^{1/2}$.
%
{
\[
K_0
\, =
\left[\ba{rrrrrrr}
\!\!  0.84 \!&\!  -0.52 \!&\!  -0.13 \!&\!  -0.07 \!&\!  -0.05 \!&\!  -0.04 \!&\!  -0.03 \!\!\\
\!\! -0.52 \!&\!   1.23 \!&\!  -0.46 \!&\!  -0.11 \!&\!  -0.06 \!&\!  -0.04 \!&\!  -0.04 \!\!\\
\!\! -0.13 \!&\!  -0.46 \!&\!   1.25 \!&\!  -0.45 \!&\!  -0.11 \!&\!  -0.06 \!&\!  -0.05 \!\!\\
\!\! -0.07 \!&\!  -0.11 \!&\!  -0.45 \!&\!   1.25 \!&\!  -0.45 \!&\!  -0.11 \!&\!  -0.07 \!\!\\
\!\! -0.05 \!&\!  -0.06 \!&\!  -0.11 \!&\!  -0.45 \!&\!   1.25 \!&\!  -0.46 \!&\!  -0.13 \!\!\\
\!\! -0.04 \!&\!  -0.04 \!&\!  -0.06 \!&\!  -0.11 \!&\!  -0.46 \!&\!   1.23 \!&\!  -0.52 \!\!\\
\!\! -0.03 \!&\!  -0.04 \!&\!  -0.05 \!&\!  -0.07 \!&\!  -0.13 \!&\!  -0.52 \!&\!   0.84 \!\!
\ea\right]
\]
\[
K_{0.01}
\, = 
\left[\ba{rrrrrrr}
\!\!  0.80 \!&\!  -0.55 \!&\!  -0.14 \!&\!   0    \!&\!   0    \!&\!   0    \!&\!  -0.11 \!\!\\
\!\! -0.55 \!&\!   1.19 \!&\!  -0.47 \!&\!  -0.17 \!&\!   0    \!&\!   0    \!&\!   0    \!\!\\
\!\! -0.14 \!&\!  -0.47 \!&\!   1.22 \!&\!  -0.45 \!&\!  -0.16 \!&\!   0    \!&\!   0    \!\!\\
\!\!  0    \!&\!  -0.17 \!&\!  -0.45 \!&\!   1.24 \!&\!  -0.45 \!&\!  -0.17 \!&\!   0    \!\!\\
\!\!  0    \!&\!   0    \!&\!  -0.16 \!&\!  -0.45 \!&\!   1.22 \!&\!  -0.47 \!&\!  -0.14 \!\!\\
\!\!  0    \!&\!   0    \!&\!   0    \!&\!  -0.17 \!&\!  -0.47 \!&\!   1.19 \!&\!  -0.55 \!\!\\
\!\! -0.11 \!&\!   0    \!&\!   0    \!&\!   0    \!&\!  -0.14 \!&\!  -0.55 \!&\!   0.80 \!\!
\ea\right]
\]
\[
K_{0.1}
\, =
\left[\ba{rrrrrrr}
\!\!  0.57 \!&\!  -0.57 \!&\!   0    \!&\!   0    \!&\!   0    \!&\!   0    \!&\!   0    \!\!\\
\!\! -0.57 \!&\!   1.14 \!&\!  -0.57 \!&\!   0    \!&\!   0    \!&\!   0    \!&\!   0    \!\!\\
\!\!  0    \!&\!  -0.57 \!&\!   1.14 \!&\!  -0.57 \!&\!   0    \!&\!   0    \!&\!   0    \!\!\\
\!\!  0    \!&\!   0    \!&\!  -0.57 \!&\!   1.14 \!&\!  -0.57 \!&\!   0    \!&\!   0    \!\!\\
\!\!  0    \!&\!   0    \!&\!   0    \!&\!  -0.57 \!&\!   1.14 \!&\!  -0.57 \!&\!   0    \!\!\\
\!\!  0    \!&\!   0    \!&\!   0    \!&\!   0    \!&\!  -0.57 \!&\!   1.14 \!&\!  -0.57 \!\!\\
\!\!  0    \!&\!   0    \!&\!   0    \!&\!   0    \!&\!   0    \!&\!  -0.57 \!&\!   0.57 \!\!
\ea\right].
\]
}
\end{example}

\section{Conclusions and Future Work}
\label{concl.sec}

We have proposed an optimization framework for the design of (sparse) interconnection graphs in LC-oscillator synchronization problems. We have identified scenarios under which the optimization problem is convex and can be solved efficiently.

Our ultimate goal is to establish a constructive framework for the synchronization of oscillator networks, in which not just the issue of synchronization but the broader questions of optimality and design of interconnection topology can be addressed. For example, it can be shown that a linearization around the consensus state of the nonlinear `swing equations,' that arise in the description of power systems, can be placed in the design framework developed in this paper and ultimately expressed as a semidefinite program. As another example, it can be shown that after applying a sequence of transformations to (\ref{CL}), the resulting equations closely resemble those of the Kuramoto oscillator. We aim to exploit these similarities for the purpose of optimal network design in our future work.

	\vspace*{-2ex}
\appendix

	\vspace*{-2ex}
\subsection*{Proof of Lemma\,\ref{zero-mode.thm}}

The proof of (i) follows from the symmetry of $\A$ and the assumption $\A \mathds{1} = 0$.

To prove (ii), note that since $\A$ is symmetric it can be diagonalized using a unitary transformation ${\cal V}$,
$
\A = {\cal V} \Lambda {\cal V}^T
$,
where $v_i$, $i = 1,\ldots,n$ denote the orthonormal eigenvectors of $\A$ and constitute the columns of ${\cal V}$; $\lambda_i$, $i = 1,\ldots,n$ denote the eigenvalues of $\A$ and constitute the diagonal elements of $\Lambda$, $\Lambda = \diag\{ \lambda_i, i = 1,\ldots,n \}$. Recalling that $\A \mathds{1} = 0$, we assume without loss of generality that
\[
\lambda_1 \,=\, 0,
~~~~~
v_1 \,=\, \mathds{1}/\sqrt{n}.
\]
Then
$
{\cal V}
=
[
\frac{1}{\sqrt{n}} \mathds{1} ~~ \widetilde{\cal V}
]
$ with $\widetilde{\cal V}^T \mathds{1} = 0$, and
$
{\cal V}^T \A {\cal V} = \diag\{ 0, \widetilde{\Lambda} \}
$,
where $\widetilde{\Lambda} = \diag\{ \lambda_i, i = 2,\ldots,n \}$. Since $\A$ is negative definite on $\mathds{1}^\perp$ then $\widetilde{\Lambda} \prec 0$. Similarly
$
{\cal V}^T \Q {\cal V} = \diag\{ 0, \widetilde{\Q} \}
$,
which results from $\Q \mathds{1} = 0$ and $\Q = \Q^T$. And since $\Q$ is positive semidefinite on $\mathds{1}^\perp$ then $\widetilde{\Q} \succeq 0$.

Multiplying $\A^T \Pee + \Pee \A = -\Q$ from the left and right by ${\cal V}^T$ and ${\cal V}$, 
and using 
\[
{\cal V}^T \Pee {\cal V}
\,=:
\left[\ba{cc}
p_1 & p_0^T \\
p_0 & \widetilde{\Pee}
\ea\right]
\]
we arrive at
\begin{align*}
\!\!
\left[\ba{cc}
0 & 0 \\
0 & \widetilde{\Lambda}
\ea\right]\!
\left[\ba{cc}
p_1 \!&\! p_0^T \! \\
p_0 \!&\! \widetilde{\Pee} \!
\ea\right]
\!+\!
\left[\ba{cc}
p_1 \!&\! p_0^T \! \\
p_0 \!&\! \widetilde{\Pee} \!
\ea\right]\!
\left[\ba{cc}
0 & 0 \\
0 & \widetilde{\Lambda}
\ea\right]\!
=
-\! \left[\ba{cc}
0 & 0 \\
0 & \widetilde{\Q}
\ea\right]\!,
\end{align*}
where $p_1$ is a scalar, $p_0$ is a column vector, and $\widetilde{\Pee}$ is a matrix. Rewriting the above equation component-wise gives
\begin{align*}
\widetilde{\Lambda} p_0
\,=\,
0,
~~~~~
\widetilde{\Lambda} \widetilde{\Pee} \,+\, \widetilde{\Pee} \widetilde{\Lambda}
\,=\,
- \widetilde{\Q},
\end{align*}
and $p_1$ is a ($\A$-- and $\Q$--independent) free parameter. From $\widetilde{\Lambda} \prec 0$ it follows that $p_0 = 0$, and therefore
$
\Pee {\cal V} = {\cal V} \, \diag\{ p_1, \widetilde{\Pee} \}
$.
In particular, this implies $\Pee \mathds{1} = p_1 \mathds{1}$ and thus $\mathds{1}$ is an eigenvector of $\Pee$ with $p_1$ as its corresponding eigenvalue. 
Furthermore, it is easy to show that
$\widetilde{\Pee} \DefinedAs \int_0^\infty e^{\widetilde{\Lambda} t} \, \widetilde{\Q} \, e^{\widetilde{\Lambda} t} \, dt$ 
is the unique solution to the Lyapunov equation 
$\widetilde{\Lambda} \widetilde{\Pee} + \widetilde{\Pee} \widetilde{\Lambda}
= - \widetilde{\Q}$
when $\widetilde{\Lambda} \prec 0$. 
Since $\widetilde{\Q} \succeq 0$ then $\widetilde{\Pee} \succeq 0$, and $\Pee$ is positive semidefinite when restricted to the subspace $\mathds{1}^\perp$.

Finally, we have
$
\trace \, (\Pee) = p_1 + \trace \, (\widetilde{\Pee})
$.
If $\Pee \succeq 0$ then $p_1 \geq 0$. 
Hence the minimum trace of $\Pee$ is achieved for $p_1 = 0$, which renders $\Pee \mathds{1} = 0$. This proves statement (ii).

To prove (iii), we note that
\begin{align*}
\Pee 
\,&=\,
{\cal V}
\left[\ba{cc}
p_1 & 0 \\
0 & \widetilde{\Pee}
\ea\right]
{\cal V}^T
\,=\,
(p_1/n) \, \mathds{1} \mathds{1}^T \,+\, \widetilde{\cal V} \widetilde{\Pee} \widetilde{\cal V}^T,
\end{align*}
with 
$\widetilde{\cal V} \widetilde{\Pee} \widetilde{\cal V}^T \mathds{1} = 0$.
From $\widetilde{\Q} \succeq 0$ and $\widetilde{\Pee} \succeq 0$ it follows that
$\Pee^\perp \DefinedAs \widetilde{\cal V} \widetilde{\Pee} \widetilde{\cal V}^T$ is positive semidefinite on $\mathds{1}^\perp$.

Also, from $\widetilde{\Lambda} \widetilde{\Pee} + \widetilde{\Pee} \widetilde{\Lambda} = - \widetilde{\Q}$ and \cite[Lemma 1]{bamdah03} we have
$\trace \, (\widetilde{\Pee}) = -(1/2) \, \trace \, (\widetilde{\Q} \widetilde{\Lambda}^{-1})$,
which implies
\begin{align*}
\trace \, (\Pee)
\,&=\,
p_1 \,+\, \trace \, (\widetilde{\Pee}) 
\,=\,
p_1 \,-\, (1/2) \, \trace \, (\widetilde{\Q} \widetilde{\Lambda}^{-1}).
\end{align*}
To express $\trace \, (\Pee)$ in terms of the Moore--Penrose pseudoinverse $\A^\dagger$, note that from $\A = {\cal V} \Lambda {\cal V}^T$ and the SVD procedure for finding the pseudoinverse, we have 
$
\A^\dagger = {\cal V} \Lambda^\dagger {\cal V}^T = {\cal V} \, \diag\{ 0, \widetilde{\Lambda}^{-1} \} {\cal V}^T
$.
Thus
$
{\cal V}^T \A^\dagger {\cal V} = \diag\{ 0, \widetilde{\Lambda}^{-1} \}
$,
and
\begin{align*}
\trace \, (\widetilde{\Q} \widetilde{\Lambda}^{-1})
\,=\,
\trace \, (\left[\ba{cc}
0 & 0 \\
0 & \widetilde{\Q} 
\ea\right]\!
\left[\ba{cc}
0 & 0 \\
0 & \widetilde{\Lambda}^{-1} \!\!\!
\ea\right]) 
\,=\,
\trace \, (\Q \A^\dagger).
\end{align*}
This proves statement (iii).

To prove (iv), 
we use \cite{ghoboysab08}
$
\A^\dagger
=
(\A - \mathds{1} \mathds{1}^T/n)^{-1} + \mathds{1} \mathds{1}^T/n,
$
which gives
\begin{align*}
\trace \, (\Q \A^\dagger)
\,&=\,
\trace \, (\Q (\A - \mathds{1} \mathds{1}^T/n)^{-1} + \Q \mathds{1} \mathds{1}^T/n)  
\,=\,
\trace \, (\Q (\A - \mathds{1} \mathds{1}^T/n)^{-1}).
\end{align*}
This proves statement (iv). The proof of the lemma is now complete.


	\vspace*{-2ex}
\subsection*{Proof of Proposition\,2}

If $Q_2$ was positive definite and thus invertible, then replacing $(1/2) \, \trace \, (Q_2^{1/2} (K + \mathds{1} \mathds{1}^T/n)^{-1} Q_2^{1/2})$ with $(1/2) \, \trace \, (X)$ in the objective function of (\ref{CVX}), subject to the linear matrix inequality (LMI) constraint
\[
\left[\ba{cc}
X & Q_2^{1/2} \\
Q_2^{1/2} & K + \mathds{1} \mathds{1}^T/n
\ea\right] \succeq 0,
\]
would follow from a simple application of the Schur complement \cite{boyelg94}. In particular, this would establish the positive definiteness, and thus the invertibility, of $K + \mathds{1} \mathds{1}^T/n$, and the optimal $X$ that minimizes the objective subject to the LMI constraint would be given by $X = Q_2^{1/2} (K + \mathds{1} \mathds{1}^T/n)^{-1} Q_2^{1/2}$. However, since $Q_2 \mathds{1} = 0$ and therefore $Q_2^{1/2}$ is singular, we can only conclude from the above LMI that $K + \mathds{1} \mathds{1}^T/n$ is positive semidefinite. We next demonstrate that $K + \mathds{1} \mathds{1}^T/n$ is indeed invertible.

Using the unitary transformation ${\cal V}$ from the proof of Lemma 1 and defining 
$
V = \diag\{ {\cal V}, {\cal V} \}
$, 
the above LMI holds if and only if
\begin{align*}
&
V^T
\! \left[\ba{cc}
X & Q_2^{1/2} \\
Q_2^{1/2} & K + \mathds{1} \mathds{1}^T/n
\ea\right]\!
V
=
\left[\ba{cccc}
x_1 & x_0^T & 0 & 0 \\
x_0 & \widetilde{X} & 0 & \widetilde{Q} \\
0 & 0 & 1 & 0 \\
0 & \widetilde{Q} & 0 & \widetilde{K}
\ea\right]
\succeq 0,
\end{align*}
where the upper-left 2-by-2 block is the appropriately partitioned matrix
$
{\cal V}^T X {\cal V}
$,
the upper-right 2-by-2 block is the partitioned matrix
$
{\cal V}^T Q_2^{1/2} {\cal V}
$
in which the zero structure follows from $Q_2^{1/2} \mathds{1} = 0$,
and the lower-right 2-by-2 block is the partitioned matrix
$
{\cal V}^T (K + \mathds{1} \mathds{1}^T/n) {\cal V}
$
in which the zero structure follows from $K \mathds{1} = 0$. Using a permutation of the rows and columns, the above LMI is equivalent to
\[
\left[\ba{cccc}
x_1 & 0 & x_0^T & 0 \\
0 & 1 & 0 & 0 \\
x_0 & 0 & \widetilde{X} & \widetilde{Q} \\
0 & 0 & \widetilde{Q} & \widetilde{K}
\ea\right]
\succeq 0,
\]
which in particular implies that
\[
\left[\ba{cc}
\widetilde{X} & \widetilde{Q} \\
\widetilde{Q} & \widetilde{K}
\ea\right]
\succeq 0.
\]
Since $Q_2^{1/2}$ is positive definite on $\mathds{1}^\perp$ then $\widetilde{Q} = \widetilde{\cal V}^T Q_2^{1/2} \widetilde{\cal V} \succ 0$. Thus $\widetilde{K} = \widetilde{\cal V}^T (K + \mathds{1} \mathds{1}^T/n) \widetilde{\cal V} \succ 0$, and therefore ${\cal V}^T (K + \mathds{1} \mathds{1}^T/n) {\cal V} \succ 0$. This implies $(K + \mathds{1} \mathds{1}^T/n) \succ 0$, and $K + \mathds{1} \mathds{1}^T/n$ is invertible.

The term $\| W \circ K \|_{\ell_1}$ in the objective function can be replaced with $\trace \, (\mathds{1} \mathds{1}^T Y)$, subject to the LMI constraint \cite{boyd}
$
-Y \leq W \circ K \leq Y
$.
Finally, the constraint $K \in \eL$, which guarantees that $K$ is the weighted Laplacian of a connected graph, is equivalent to the set of conditions
$K = K^T$, $M \circ K \leq 0$, $K \mathds{1} = 0$, $K + \mathds{1} \mathds{1}^T/n \succ 0$.
The first and last of these conditions are automatically fulfilled when $K$ satisfies the LMI, and are therefore dropped from the formulation of the optimization problem. The proof of the proposition is now complete.


\end{document}